\documentclass[12pt]{article}
\usepackage{latexsym,subfigure}
\usepackage{comment,tikz}
\usepackage{showlabels}
\usepackage{amsmath,amsthm,amsfonts,amssymb,graphicx,epsfig,latexsym,color}
\usepackage{epsf,enumerate,url}

\newtheorem{thm}{Theorem}

{}
{}

\theoremstyle{remark}

\newtheorem{remark}[thm]{Remark}

\newtheorem*{definition*}{Definition}
\newtheorem*{remark*}{Remark}

\def\R{{\mathbb R}}

\def\Z{{\mathbb Z}}

\begin{document}

\author{Koji Fujiwara
\thanks{The author is supported in part by
Grant-in-Aid for Scientific Research (No.15H05739, 20H00114). }
}

\title{An example of a closed 5-manifold of nonpositive curvature that fibers over a circle}

\maketitle 

\begin{abstract}
We exhibit a closed aspherical 5-manifold of nonpositive curvature that fibers over a circle whose fundamental group is hyperbolic relative to abelian subgroups such that the fiber is a closed aspherical 4-manifold whose fundamental group is not hyperbolic relative to abelian subgroups. 
\end{abstract}

Recently interesting examples of hyperbolic groups with subgroups of finite type
that are not hyperbolic are constructed by Italiano, Martelli, and Migliorini \cite[Corollary 2]{IMM}.
They start with a hyperbolic 5-manifold
of finite volume with cusps that fibers over a circle, \cite[Theorem 1]{IMM}. Then they close
the cusps using {\it $2\pi$-fillings} given in \cite{FM} and 
obtain locally CAT($-k$), $k>0$, closed pseudo-manifolds that fibers over a circle,
whose fundamental group is hyperbolic, such that the fundamental
group of the fiber is not hyperbolic. 

In this note, we use different $2\pi$-fillings to close the cusps of their
cusped manifold, which gives a closed 5-manifold of nonpositive curvature
that fibers over a circle such that the fiber is a closed apherical manifold.

\begin{thm}
  
There exists a 5-dimensional closed Riemannian manifold $V$
of nonpositive sectional curvature that fibers over $S^1$ 
such that 
\begin{enumerate}
\item
$V$ has disjoint totally geomedesic flat $3$-tori such that the sectional curvature is negative on all 2-planes that are not tanget to those tori. 
The universal cover $\tilde V$ is an Hadamard manifold with disjoint colletion 
of 3-dimensional Euclidean spaces, $\Bbb E^3$ (called ``isolated flats''). 

\item
$G=\pi_1(V)$ is hyperbolic relative to subgroups
isomorphic to $\Z^3$, which are the fundamental groups of 
the $3$-tori. 
\item
The fiber, $W$, is a 4-dimensional, closed, aspherical manifold and 
$H=\pi_1(W)$ is not hyperbolic relative to any abelian subgroups. 
\end{enumerate}

\end{thm}

\proof
{\it Construction}.
The construction is very similar to the one 
for \cite[Corollary 2]{IMM}. We review their construction. 
They start with the hyperbolic 5-manifold $M$ with cusps that fibers
over $S^1$ with the (connected) fiber $F$ given in \cite[Section 1]{IMM}.
In the section 3, truncating at every cusp of $M$ they 
obtain  a compact manifold $\bar M$ whose boundary consists of flat 4-tori.
$\bar M$ fibers over a circle with a fiber $\bar F=F \cap \bar M$.
Replacing $M$ with a finite cover if necessary, they
arrange that each boundary flat 4-torus satisfies the {\it $2\pi$-condition}, namely, any non-trivial (in the fundamental group)
 curve on the tori
has length larger than $2\pi$. 
They also arrange that on each boundary torus, the fibration of $M$ restricts to a fibration 
over a circle such that each fiber is a geodesic flat 3-torus. 

Then, gluing a {\it $2\pi$-filling}
in the sense of \cite{FM} to each boundary 4-torus of $\bar M$, 
they obtain a pseudo-manifold
$\hat M$ that fibers over $S^1$ with the fiber $\hat F$.
Topologically, on each boundary 4-torus, 
they cone off (shrink) each 3-torus fiber to a point, where $\hat F$
is obtained from $\bar F$ by shrinking each boundary 3-torus to a point. 
By the terminology of $\cite{FM}$, their $2\pi$-filling is 
the {\it partial cone} $C(T^4, T^3)$.
By \cite[Theorem 2.7]{FM}, $\hat M$ is locally CAT($-k$) with $k>0$, 
therefore, $\pi_1(\hat M)$ is a hyperbolic group. The subgroup $\pi_1(\hat F)$ 
is of finite type since $\hat F$ is aspherical, and they argue that it is not hyperbolic.  
Both $\hat M$ and $\hat F$ are orientable.

Our manifold $V$ is given by a variant of their construction. 
We start with their manifold $M$.
To obtain a manifold, we use a  different $2\pi$-filling to 
each boundary torus of $\bar M$.
Instead of $C(T^4,T^3)$, we glue
the partial cone $C(T^4,S^1)$ to each boundary $T^4$ of $\bar M$, 
where its {\it core} $V(T^4,S^1)$ is a flat 3-torus. 
In their case, the core $V(T^4,T^3)$ is a circle, where the pseudo-manifold $\hat M$
has singularity. 

We recall the construction of $C(T^4,S^1)$.
Let $T^4$ be a flat torus in the boundary of $\bar M$,
which fibers over a circle with flat 3-tori fibers. 
We further foliate each 3-torus fiber $T^3$ by simple
closed geodesics, parellel to one, $\gamma$,  that represents an element, $g \in \pi_1(T^3)\simeq \Bbb Z ^3$,
with $\Bbb Z^3/\langle g \rangle =\Bbb Z^2$.
Topologically, we then cone off each of those loops to a point,
and obtain $C(T^4,S^1)$.
The collection of those cone points, the core,  will give a flat 3-torus in $V$.

$V$ still fibers over $S^1$ with the fiber $W$, which 
is obtained from $F\cap \bar M$ by closing each boundary 3-torus to a 2-torus. $W$ is a closed $4$-manifold. 
Both $V$ and $W$ are orientable.

By the {\it cusp closing theorem} by Schroeder, \cite{S},
the closed manifold $V$ has a Riemannian metric of nonpositive sectional curvature such that it has a totally geodesic flat 3-torus from each cusp of $M$, which gives a disjoint collection of tori, and that 
the sectional curvature is strictly negative on all 2-planes that are not tangent to those tori. 
One could  view this as a special case of \cite[Theorem 2.7]{FM}, 
see Remark 2.9 in there, although they do not explicitly say that 
the metric is Riemannian.

$\tilde V$ is an Hadamard manifold with {\it isolated flats} in the sense 
of Hruska-Kleiner, \cite{HK} that are lifts
of the 3-tori in $M$.
$V$ is aspherical, and 
since $V$ is covered by $W \times \R$, $W$ is also aspherical.

Let $T_i$ be the embedded 3-tori in $V$ and set $P_i=\pi_1(T_i)$.
$G=\pi_1(V)$ is hyperbolic relative to $\{P_i\}$, \cite{HK}. 
Such group is called {\it toral relatively hyperbolic}
since the subgroups $P_i$ are abelian. $G$ is a torsion free, CAT(0) group.

{\it The fiber group}.
Let $H=\pi_1(W)$.
It remains to argue that 
$H$ is not hyperbolic relative to any abelian subgroups.
We refer to \cite{H} for the definitions of relative hyperbolicity.
First, any non-trivial abelian subgroups in $G$ is $\Z,\Z^2$, or $\Z^3$ by the construction of $V$
and the flat torus theorem by Lawson-Yau.
It implies that any non-trivial abelian subgroup in $H$ is $\Z$ or $ \Z^2$, since 
if $H$ has a subgroup isomorphic to $\Z^3$, then it must be conjugate to a subgroup of finite index in one of $\pi_1(T_i)$,  then 
it has a non-trivial projection to $\pi_1(S^1)$, impossible.

To argue by contradiction, 
suppose $H$ is hyperbolic relative to a collection of abelian subgroups $\{A_i\}$.
We may remove all $A_i$ that is isomorphic to $\Z$
without destroying relative hyperbolicity (see \cite{GL}), so that we assume each $A_i$
is $\Z^2$.
Let $Q_i=H \cap P_i$, which is $\Z^2$. 

Each $Q_i$ is a maximal abelian subgroup
in $H$. This is because if $R_i$ is abelian and $Q_i<R_i<H$, then 
$R_i$ must be conjugate into one of $P_j$ in $G$, by the flat torus
theorem, which has to be $P_i$. But since $Q_i=H \cap P_i$,
$Q_i=R_i$.

Each $Q_i$ is 
conjugate to one of $\{A_i\}$  in $H$. 
This is because
since each $Q_i$ is $\Z^2$, it must be {\it parabolic}
with respect to $(H, \{A_i\})$, which means $Q_i$ is a subgroup
of some $A_j$ after a conjugation in $H$.
Since $Q_i$ is maximal abelian in $H$, it is in fact  conjugate to $A_j$.

Let $b\in G$ be an element that maps to a generator of $\pi_1(S^1)$.
By conjugation, $b$ defines an element, denoted by $b$, in ${\rm Aut}(H)$.
It will map each $A_i$ to a conjugate of some $A_j$ in $H$, so that there exists $L>0$ s.t.
$b^L$ is in ${\rm Aut}(H, \{A_i\})$, whose elements
preserves each $A_i$ upto conjugation in $H$
(see \cite{GL}).
We claim $b^L$ has infinite order in ${\rm Out}(H,\{A_i\})$.
Indeed, suppose not. Then for some $k\not=0$, there exists $a\in H$
such that for all $h\in H$, $b^khb^{-k}=aha^{-1}$.
Namely, $a^{-1}b^k$ centralizes $H$.
Now since each $A_i$ is malnormal in $H$, 
$a^{-1} b^k\in \cap_i Q_i$.
But $\cap_i Q_i$ is trivial, which means
$b^k =a\in H$, impossible since $b^k \not\in H$.

Since $b^L$ has infinite order, ${\rm Out}(H,\{A_i\})$
is infinite. 
Now by \cite[Theorem 1.3]{GL}, $H$ must split over a cyclic or a parabolic subgroup, which 
is $\Z$ or $\Z^2$.
As in \cite{IMM}, an easy Mayer-Vietoris type argument implies $H^4(W)=0$,
which is a contradiction since $W$ is closed. 
\qed

We conclude with some remarks. 
\begin{remark}
\begin{enumerate}
\item
The fundamental group of the fiber $\bar F$ of 
a pseudo-manifold $\bar M$  in \cite{IMM} is not toral relatively 
hyperbolic either (since it is not hyperbolic and does not contain non-trivial 
abelian subgroups).

The fundamental group of their cusped manifold $M$ is toral relatively hyperbolic with 
respect to the cusp subgroups $\Bbb Z^4$.
Jason Manning pointed out to the author that one could  show using \cite{MW} that the fundamental 
group of the fiber $F$ is not toral relatively hyperbolic. A new aspect of our example is that both $G$ and $H$ are the fundamental groups of 
closed aspherical manifolds.

\item In \cite{IMM}, they take a finite cover of $M$ to make sure that 
$\bar M$ satisfies the $2\pi$-condition. But for our construction,
we do not have to take a finite cover, since one can always find
a simple closed geodesic $\gamma$ in the fiber $3$-torus in each boundary  $4$-torus
of $\bar M$, which is enough for our argument. 
\end{enumerate}
\end{remark}

{\bf Acknowledgement}
The author would like to thank Danny Calegari for bringing the work \cite{IMM}
to his attention. 
He is grateful to Vincent Guirardel, Jason Manning and  Mladen Bestvina
for useful comments.

\end{document}